\documentclass[12pt]{amsart}
\usepackage{amssymb}
\usepackage{times}

\theoremstyle{plain}

\numberwithin{equation}{section}

\newcommand{\calA}{\mathcal{A}}
\newcommand{\calB}{\mathcal{B}}
\newcommand{\calC}{\mathcal{C}}
\newcommand{\calD}{\mathcal{D}}

\newcommand{\calL}{\mathcal{L}}

\newcommand{\calO}{\mathcal{O}}
\newcommand{\calP}{\mathcal{P}}

\newcommand{\bbF}{\mathbb{F}}

\newcommand{\bbP}{\mathbb{P}}
\newcommand{\bbQ}{\mathbb{Q}}
\newcommand{\bbR}{\mathbb{R}}

\newcommand{\bbZ}{\mathbb{Z}}

\def\Aut{\text{Aut}}

\def\Pic{\text{Pic}}

\def\PSL{\text{PSL}}
\def\PGL{\text{PGL}}

\begin{document}
\title [A supersingular K3 surface and the Leech lattice] {A supersingular
K3 surface in
characteristic 2 and the Leech lattice}
\author{I. Dolgachev}
\address{Department of Mathematics, University of Michigan, Ann Arbor, MI
48109,USA}
\email{idolga@umich.edu}
\thanks{Research of the first author was partially supported by NSF grant
DMS 9970460 and the Japanese Society
for Promotion of Science}

\author{S. Kond{$\bar{\rm o}$}}
\address{Graduate School of Mathematics, Nagoya University, Nagoya,
464-8602, Japan}
\email{kondo@math.nagoya-u.ac.jp}
\thanks{Research of the second author is partially supported by
Grant-in-Aid for Exploratory Research 11874004} 

\begin{abstract}
We construct a K3 surface over an algebraically closed field
of characteristic 2 which contains two sets of 21 disjoint smooth rational
curves such that each curve from one set intersects exactly 5 curves from
the other set. This configuration is isomorphic to the configuration of
points and lines on the projective plane over the finite field of 4
elements. The surface admits a finite automorphism group isomorphic to
$\PGL(3,\bbF_4)\cdot 2$ such that a subgroup $\PGL(3,\bbF_4)$ acts on the
configuration of each set of
21 smooth rational curves, and the additional element of order 2
interchanges the two sets. The Picard lattice of the surface
is a reflective sublattice of an even unimodular lattice $II_{1,25}$ of
signatuire
$(1,25)$ and the classes of the 42 curves correspond to some Leech roots in
$II_{1,25}$

\end{abstract}
\maketitle

\section{Introduction}
Let $k$ be an algebraically closed field of characteristic 2. Consider
$\bbF_4\subset k$ and
$\bbP^2(\bbF_4)\subset \bbP^2(k)$. Let $\calP$ be the set of points and let
$\check{\calP}$ be the set of lines in $\bbP^2(\bbF_4)$. Each set contains
21 elements, each point is contained in exactly 5 lines and each line
contains exactly 5 points. It is known that the group of automorphisms of
the configuration $(\calP,\check{\calP})$ is isomorphic to $M_{21} \cdot
D_{12}$ where
$M_{21} ( \cong \PSL(3,\bbF_4) )$ is a simple subgroup of the Mathieu group
$M_{24}$ and $D_{12}$ is a dihedral group of order 12. 

In this paper we prove the following theorem: 

\subsection{Main Theorem}\label{main} {\it There exists a unique (up to
isomorphism) $K3$ surface over
$k$ satisfying the following equivalent properties $:$ \begin{itemize}
\item[(i)]  The Picard lattice of $X$ is isomorphic to $U\perp
D_{20}
\ ;$
\item[(ii)] $X$ has a jacobian quasi-elliptic fibration with one fibre of
type $\tilde{D}_{20} \ ;$
\item[(iii)] $X$ has a quasi-elliptic fibration with the Weierstrass
equation \[y^2 = x^3+t^{11};\]
\item[(iv)] $X$ has a quasi-elliptic fibration with 5 fibres of type
$\tilde{D}_4$ and the group of sections isomorphic to $(\bbZ/2)^4\ ;$
\item[(v)] $X$ contains a set $\calA$ of \ $21$ disjoint $(-2)$-curves and
another set $\calB$ of \ $21$ disjoint
$(-2)$-curves such
that each curve from one set intersects exactly 5 curves from the other set
with multiplicity $1 ;$
\item[(vi)] $X$ is birationally isomorphic to the inseparable double cover
of \ $\bbP^2$ with branch divisor
\[x_0x_1x_2(x_0^3+x_1^3+x_2^3) = 0;\]
\item[(vii)] $X$ is isomorphic to a minimal nonsingular model of the
quartic surface with $7$ rational double points of type $A_3$ which is
defined by the equation
\[x_0^4+x_1^4+x_2^4+x_3^4+x_0^2x_1^2+x_0^2x_2^2+x_1^2x_2^2+x_0x_1x_2(x_0+x_1
+x_2) = 0;\]
\item[(viii)] $X$ is isomorphic to the surface in $\bbP^2\times \bbP^2$
given by the equations
\[x_0y_0^2+x_1y_1^2+x_2y_2^2 = 0,\quad x_0^2y_0+x_1^2y_1+x_2^2y_2 = 0.\]
\end{itemize}
The automorphism group $\Aut(X)$ contains a normal infinite subgroup
generated by $168$ involutions and the quotient is a finite group
isomorphic to $\PGL(3,\bbF_4)\cdot 2$.}

\medskip
Since the group $\PSL(3,\bbF_4)$ is not a subgroup of the Mathieu group
$M_{23}$, our theorem shows that the classification of finite groups acting
symplectically on a complex K3 surface due to S. Mukai \cite{Mu} does not
extend to positive characteristic.

\section{Supersingular K3 surfaces in characteristic 2}\label{s1}
\subsection{Known facts}\label{known} Recall that a {\it supersingular $K3$
surface} (in the sense of Shioda)
is a K3 surface with the Picard group of rank 22. This occurs only if the
characteristic of the ground field $k$ is a positive prime $p$. By a result
of M. Artin \cite{Artin}, the Picard lattice $S_X = \Pic(X)$ of a
supersingular K3 surface is a $p$-elementary lattice (i.e. the discriminant
group $S_X^*/S_X$ is a
$p$-elementary abelian group). The dimension $r$ of the discriminant group
over $\bbF_p$ is even and the number
$\sigma = r/2$ is called the {\it Artin invariant}. We shall assume that $p
= 2$. A fundamental theorem of Rudakov and Shafarevich
\cite{RS1},\cite{RS2} tells that any supersingular K3 surface admits a
quasi-elliptic fibration, i.e. a morphism $f:X\to\bbP^1$ whose general
fibre is a regular but not smooth geometrically irreducible curve of genus
1. Over an open subset $U$ of the base, each fibre is isomorphic to an
irreducible cuspidal cubic, the reducible fibres are Kodaira genus 1 curves
of additive type. The closure of the set of cusps of irreducible fibres is
a smooth irreducible curve $C$, the {\it cusp curve}. The restriction of
$f$ to $C$ is a purely inseparable cover of degree 2. It follows that $C$
is isomorphic to $\bbP^1$, hence by adjunction is a $(-2)$-curve on $X$. A
surface with a quasi-elliptic fibration is unirational. Thus any
supersingular K3 surface in characteristic 2 is unirational. This is one of
the main results of
\cite{RS1},\cite{RS2}.

In this paper we study a supersingular K3 surface with the Artin invariant
$\sigma = 1$. It follows from the classification of 2-elementary lattices
of signature $(1,21)$ that the Picard lattice of such a K3 surface is
unique up to isometries (see \cite{RS2}, \S 1). Therefore it is isomorphic
to $U \perp D_{20}$. Here, as usual, $U$ denotes the unique even unimodular
indefinite lattice of rank 2 and $A_n,D_n,E_n$ denote the negative definite
even lattices defined by the opposite of the Cartan matrix of the simple
root system of the corresponding type. Moreover it is known that any
supersingular K3 surface with Artin invariant $\sigma = 1$ is unique up to
isomorphisms (\cite{RS2}, \S 11).

\subsection{The Weierstrass equation}\label{lattice} 

\subsection{Proposition}\label{P1} {\it Let $X$ be a supersingular $K3$
surface whose Picard lattice $S_X$ is isometric to $U\perp D_{20}$. Then
$X$ has a quasi-elliptic fibration with one singular fibre of type
$\tilde{D}_{20}$ and a section. Its Weierstrass equation is
\[y^2+x^3+t^{11} = 0.\]}

\begin{proof} First we shall see that the above Weierstrass equation
defines a K3 surface whose Picard lattice is isomorphic to $U\perp D_{20}$.
The reducible fibres lie over points at which the affine surface has a
singular point.
Taking the derivatives we find that such point is given by $x = t = 0$.
This gives $t = 0$ as the only solution for ''finite'' $t$. When $t =
\infty$, we make the substitutions
$y\to t^{6}y, x\to t^4x$ and then $t = 1/\tau$. Then we have the equation
\[y^2+x^3+\tau = 0 \] which is nonsingular over $\tau = 0$. The
discriminant $\Delta$ of the fibration has order 20 at $t = 0$. Hence the
fibration defines a K3 surface
(see \cite{RS2}, \S 12, \cite{CD}, Chapter 5). It has a unique reducible
fibre. It follows from the classification of fibres of quasi-elliptic
fibrations that it must be of type $\tilde{D}_{20}$. Obviously its Picard
lattice contains $U\perp D_{20}$ which is generated by the class of
components of fibres and a section. Since $\sigma \geq 1$, the Picard
lattice isomorphic
to $U\perp D_{20}$.
Now the assertion follows from the uniqueness of supersingular K3 surface
with the Artin invariant 1.
\end{proof}

\section{Leech roots}\label{leech}

In this section, we denote by $X$ the K3 surface whose Picard lattice is
isomorphic to $U\perp D_{20}$.

\subsection{The Leech lattice} We follow the notation and the main ideas
from \cite{Borcherds}, \cite{Conway1},
\cite{Conway2}, \cite{Kondo}. First we embed the Picard lattice $S_X\cong
U\perp D_{20}$ in the lattice $L = \Lambda\perp U \cong II_{1,25}$, where
$\Lambda$ is the
{\it Leech lattice} and $U$ is the hyperbolic plane. We denote each vector
$x\in L$ by $(\lambda,m,n)$ where
$\lambda\in \Lambda$, and $x = \lambda+mf+ng,$ with $f,g$ being the
standard generators of $U$, i.e. $f^2 = g^2 = 0$, and $\langle f, g \rangle
= 1$.
Note that $r = (\lambda,1,-1-\frac{\langle \lambda,\lambda \rangle} {2})$
satisfies $r^{2} = -2$.
Such vectors will be called {\it Leech roots}. We denote by $\Delta(L)$ the
set of all Leech roots.
Recall that $\Lambda$ can be defined as a certain lattice in $\bbR^{24} =
\bbR^{\bbP^1(\bbF_{23})}$ equipped with inner product $\langle x,y \rangle
= -\frac{x\cdot y}{8}$. For any subset $A$ of $\Omega = \bbP^1(\bbF_{23})$
let $\nu_A$ denote the vector $\sum_{i\in A} e_i$, where $\{e_\infty,
e_0,\ldots,e_{22}\}$ is the standard basis in $\bbR^{24}$. A {\it Steiner
system} $S(5,8,24)$ is a set consisting of eight-element subsets of
$\Omega$ such that
any five-element subset belongs to a unique element of $S(5,8,24)$. An
eight-element subset in $S(5,8,24)$ is called an {\it octad}. Then
$\Lambda$ is
defined as a lattice generated by the vectors $\nu_\Omega-4\nu_\infty$ and
$2\nu_K$, where $K$ belongs to the
Steiner system $S(5,8,24)$. Let $W(L)$ be the subgroup generated by
reflections in the orthogonal group $O(L)$ of $L$. Let $P(L)$ be a
connected component of
$$P(L) = \{ x \in \bbP(L\otimes \bbR) : \langle x,x \rangle > 0 \}.$$ Then
$W(L)$ acts naturally on $P(L)$. A fundamental domain of this action of
$W(L)$ is given by
$$\calD = \{ x \in P(L) : \langle x, r \rangle > 0, \ r\in \Delta(L) \}.$$
It is known that $O(L)$ is a split extension of $W(L)$ by $\Aut(\calD)$
(\cite{Conway2}).

\subsection{Lemma}\label{l1} {\it Let $X$ be the $K3$ surface whose Picard
lattice $S_X$ is isomorphic to $U\perp D_{20}$. Then there is a primitive
embedding of $S_X$ in $L$ such that the orthogonal complement $S_X^\perp$
is generated by some Leech roots and isomorphic to the root lattice $D_4$.}

\begin{proof} Consider the following vectors in $\Lambda$:
\begin{equation}\label{XYZ}
X = 4\nu_{\infty} + \nu_{\Omega},\quad Y = 4\nu_{0} + \nu_{\Omega},\quad Z
= 0,\end{equation}
\begin{equation}
T = (x_{\infty}, x_{0}, x_{1}, x_{k_{2}},..., x_{k_{22}}) =
(3,3,3,-1,-1,-1,-1,-1,1,...,1)\end{equation} where $K = \{ \infty, 0,
1,k_{2},...,k_{6} \}$ is an octad. The corresponding Leech roots
\begin{equation}\label{xyz}
x = (X, 1, 2), \quad y = (Y, 1, 2),\quad z = (0, 1, -1),\quad t = (T, 1, 2)
\end{equation}
generate a root lattice $R$ isomorphic to $D_{4}.$ Obviously $R$ is
primitive in $L$.

Let $S$ be the orthogonal complement of $R$ in $L$. Then it is an even
lattice of signature $(1,21)$ and the discriminant group isomorphic to
$(\bbZ/2\bbZ)^2$ (the discriminant groups of $S$ and $R$ are isomorphic).
Since such a lattice is unique, up to isometry, we obtain that $S \cong U
\perp D_{20}$.
\end{proof}

\subsection{(42 + 168) Leech roots}\label{domain} We fix an embedding of
$S_X$ into $L$ as in Lemma \ref{l1}. Let $P(X)$ be the positive cone of
$X$, i.e. a connected component of $$\{ x \in \bbP(S_X\otimes \bbR) :
\langle x,x \rangle > 0 \}$$ which contains the class of an ample divisor.
Let $\calD (X)$ be the intersection of $\calD$ and $P(X)$. It is known that
$\Aut(\calD (X)) \cong M_{21} \cdot D_{12}$ where $M_{21} \cong
\PSL(3,\bbF_4)$ and $D_{12}$ is a
dihedral group of order 12 (\cite{Borcherds}, \S 8, Example 5). Each
hyperplane bounding $\calD (X)$ is the one perpendicular to a negative norm
vector in $S_X$ which is the projection of a Leech root. There are two
possibilities of such Leech roots $r$:

\begin{itemize}
\item[(i)] $r$ and $R$ ($\cong D_4)$ generate a root lattice isomorphic to
$D_4\perp A_1$;
\item[(ii)] $r$ and $R$ generate a root lattice isomorphic to $D_5$.
\end{itemize}

\subsection{Lemma}\label{l2} {\it There are exactly forty-two Leech roots
which are orthogonal to $R$.}

\begin{proof} Let $R'$ be the sublattice of $L$ generated by vectors $x, y,
z$ from \eqref{xyz}. Obviously $R'$ is
isomorphic to the root lattice $A_3$. By Conway \cite{Conway1}, the Leech
roots orthogonal to $R'$ are
\[4\nu_{\infty} + 4\nu_{0}, \quad \nu_{\Omega} - 4\nu_{k}, \quad
2\nu_{K'},\] where
$K'$ contains $\{\infty, 0\}$ and $k \in \{1,2,...,22 \}$.

Among these Leech roots, the followings are orthogonal to vector $t\in R$:
\[4\nu_{\infty} + 4\nu_{0}, \quad \nu_{\Omega} - 4\nu_{k},\quad
2\nu_{K'},\] where $K' \in S(5,8,24)$, $\ K' \cap K = \{\infty, 0 \}$ or
$K' \cap K = \{\infty, 0, 1, * \}$ and $k \in K \setminus \{\infty, 0, 1 \}
$. As is easy to see the number of these roots is equal to 42. \end{proof} 

\subsection{Remark}\label{rk1} In his paper \cite{Todd}, Todd listed the
$759$ octads of the Steiner system $S(5,8,24)$ (p.219, Table I). The octads
in the proof of Lemma \ref{l2} correspond to the following $\{ E_{\alpha},
L_{\beta} \}$ in his table. Here we assume \[K = \{ \infty, 0, 1, 2, 3, 5,
14, 17 \}.\] Then \[E_{1} = \{ \infty, 0, 1, 2, 4, 13, 16, 22 \}, \quad
E_{2} = \{ \infty, 0, 1, 2, 6,7,19,21 \}, \] \[E_{3} = \{ \infty, 0, 1,
2,8,11,12,18 \},\quad E_{4} = \{ \infty, 0, 1, 2,9,10,15,20 \},\] \[E_{5} =
\{ \infty, 0, 1, 3,4,11,19,20 \},\quad E_{6} = \{ \infty, 0, 1, 3,6,8,10,13
\},\] \[ E_{7} = \{ \infty, 0, 1, 3,7,9,16,18 \}, \quad E_{8} = \{ \infty,
0,1,3,12,15,21,22 \}, \] \[E_{9} = \{ \infty, 0, 1,4,5,7,8,15 \}, \quad
E_{10} = \{ \infty, 0, 1, 4,6,9,12,17 \},\] \[E_{11} = \{ \infty, 0,
1,4,10,14,18,21 \},\quad E_{12} = \{ \infty, 0, 1, 5,6,18,20,22 \},\]
\[E_{13} = \{ \infty, 0, 1, 5,9,11,13,21 \}, \quad E_{14} = \{ \infty, 0,
1, 5,10,12,16,19 \},\] \[E_{15} = \{ \infty, 0, 1, 6,11,14,15,16 \},\quad
E_{16} = \{ \infty, 0, 1, 7,10,11,17,22 \}, \] \[E_{17} = \{ \infty, 0, 1,
7, 12,13,14,20 \},\quad E_{18} = \{ \infty, 0, 1, 8, 9, 14, 19, 22 \},\]
\[E_{19} = \{ \infty, 0, 1, 8, 16, 17, 20, 21 \}, \quad E_{20} = \{ \infty,
0, 1, 13, 15, 17, 18, 19 \}, \]
\[L_{1} = \{ \infty, 0, 4, 6,8,16,18,19 \}, \quad L_{2} = \{ \infty, 0,
4,6,13,15,20,21 \},\]
\[L_{3} = \{ \infty, 0, 4,7,9,10,13,19 \},\quad L_{4} = \{ \infty, 0,
4,7,11,12,16,21 \},\]
\[L_{5} = \{ \infty, 0, 4,8,10,12,20,22 \}, \quad L_{6} = \{ \infty, 0,
4,9,11,15,18,22 \}, \]
\[L_{7} = \{ \infty, 0, 6,7,8,9,11,20 \},\quad L_{8} = \{ \infty, 0,
6,7,10,12,15,18 \},\]
\[L_{9} = \{ \infty, 0,6,9,10,16,21,22 \},\quad L_{10} = \{ \infty, 0,
6,9,10,16,21,22 \}, \]
\[L_{11} = \{ \infty, 0, 7,8,13,18,21,22 \},\quad L_{12} = \{ \infty, 0,
7,15,16,19,20,22 \}, \]
\[L_{13} = \{ \infty, 0,8,9,12,13,15,16 \},\quad L_{14} = \{ \infty, 0,
8,10,11,15,19,21 \},\]
\[L_{15} = \{ \infty, 0, 9,12,18,19,20,21 \}, \quad L_{16} = \{ \infty, 0,
10,11,13,16,18,20 \}.\]
The remaining 6 Leech roots correspond to \[4\nu_{\infty} + 4\nu_{0},\quad
\nu_{\Omega} - 4\nu_{k},\]
where $k \in \{ 2,3,5,14,17 \}$.

\subsection{Lemma}\label{168} {\it There are exactly $168$ Leech roots $r$
such that $r$ and $R$ generate
a root lattice $D_5$.}

\begin{proof} We count the number of such Leech roots $r$ with $\langle r,
t \rangle = 1$.
Such $r$ corresponds to one of the following vectors in $\Lambda$ (see the
proof of Lemma \ref{l2}):
$$\nu_{\Omega} - 4\nu_{k} \ (k \notin K), \quad \nu_{K'} \ (K' \in
S(5,8,24), \ \mid K \cap K' \mid = 4,\ 1 \notin K').$$ Obviously the number
of vectors of the first type is 16. Since the number of octads containing
fixed 4 points is 5,
the number of vectors of the second type is 40. Thus we have the desired
number $56 \times 3 = 168$.
\end{proof}

\medskip
\subsection{The 42 smooth rational curves}\label{42} Let \[w = (0,0,1)\] be
the {\it Weyl vector} in $L$ (characterized by the property that $\langle
w,\lambda \rangle = 1$
for each Leech root). Since the Leech lattice does not contain
$(-2)$-vectors, $\langle w, r \rangle \not= 0$ for any $(-2)$-vector $r$ in
$L = \Lambda \perp U$. Consider its orthogonal projection $w'$ to
$R^\perp\otimes \bbQ$. Easy computation gives \[w' = w+5z+3x+3y+3t\in
R^\perp,\]
and
\begin{equation}\label{w'}
\langle w',w' \rangle = 14.
\end{equation}
Fix an isometry from $S_X$ to $R^\perp$ and let $h$ be the divisor class
corresponding to $w'$. The above property of $w$ implies that $\langle h, r
\rangle \not= 0$ for any $(-2)$-vectors $r$ in $S_X$. Composing the
embedding with
reflections in $(-2)$-vectors in $S_X$, we may assume that $h$ is an ample
divisor class. Each of the 42 Leech vectors from Lemma \ref{l2} defines a
vector $v$ from $\Pic(X)$ with self-intersection $-2$. Since $\langle h, v
\rangle = 1$, by Riemann-Roch, we obtain that $v$ is the divisor class of a
curve $R_v$ with $R_v^2 = -2$. Since $h$ is ample, each $R_v$ is an
irreducible curve, and hence isomorphic to $\bbP^1$. 

Thus the 42 Leech roots in Lemma \ref{l2} define 42 smooth rational curves
in $X$.

\subsection{Lemma}\label{l3}{\it Let $X$ be a K3 surface over a field of
characteristic 2. The following properties are equaivalent: \begin{itemize}
\item [(i)] The Picard lattice $S_X$ of
$X$ is isomorphic to $U \perp D_{20}$;
\item [(ii)]
$X$ has a quasi-elliptic fibration with five singular fibres of type
$\tilde{D}_{4}$ and sixteen disjoint sections; \item [(iii)] There are two
families $\calA$ and $\calB$ each consisting of $21$ disjoint
smooth rational curves. Each member in one family meets exactly five
members in another family. The set $\calA\cup\calB$ generates $S_X$.
\end{itemize}}

\begin{proof} (i) $\Rightarrow$ (ii) Fix an isometry from $S_X$ to
$R^\perp$. Consider the five disjoint $(-2)$-curves $K_{k}$ corresponding
to
$$\nu_{\Omega} - 4\nu_{k}, \quad ( k \in \{ 2,3,5,14,17 \}).$$ Then $\mid
2K_{2} + R_{1} + R_{2} + R_{3} + R_{4} \mid$ gives a genus 1 fibration with
five singular fibers of type $\tilde{D}_{4}$, where $R_{i}$ correspond to
the vectors $E_i$ in Remark \ref{rk1}. Moreover the curves corresponding to
the vectors $L_{i}$ are 16 disjoint sections of this fibration. The
fibration satisfies the Rudakov-Shafarevich criterion of quasi-elliptic
fibration \cite{RS2},\S 4. 

(ii) $\Rightarrow$ (i) It follows from the Shioda-Tate formula
\cite{Shioda1} that the discriminant of $S_X$ is equal to 2. 

(ii) $\Rightarrow$ (iii) We take for $\calA$ the set of sections and
multiple components of fibres. We take for $\calB$ the set of non-multiple
components of reducible fibres and the cusp curve. The curves from the set
$\calA$ correspond to the Leech roots $\{ E_{\alpha}, 4\nu_{\infty} +
4\nu_{0}
\},$ and the curves from $\calB$ correspond to the Leech roots $\{
L_{\alpha}, \nu_{\Omega} - 4\nu_{k}
\}$. The cusp curve $C$ correspond to the Leech root $4\nu_{\infty} +
4\nu_0$ (see the following Remark \ref{rk2}).
It follows from the Shioda-Tate formula that the set $\calA\cup\calB$
generates $S_X$.

(iii) $\Rightarrow$ (ii) Take
$R_0\in
\calA$ and 4 curves
$R_1,R_2,R_3,R_4$ from $\calB$ which intersect $R_0$. Then
\begin{equation}\label{fiber}
F = 2R_0+R_1+R_2+R_3+R_4
\end{equation}
defines a quasi-elliptic fibration with a fibre $F$ of type $\tilde{D}_4$.
Let $N$ be the fifth curve from $\calB$ which intersects $R_0$. It is easy
to see that the curves from $\calA\cup \calB\setminus \{N\}$ which do not
intersect the curves $R_1,R_2,R_3,R_4$ form four more reducible fibres of
type $\tilde{D}_4$.  The remaining curves give 16 disjoint sections.

\end{proof}

\subsection{Lemma}\label{l4}{\it Let $X$ be the $K3$ surface whose Picard
lattice is isomorphic to $U \perp D_{20}$ and $\calA, \calB$ are the sets
described in the previous lemma. Let $h$ be the divisor class corresponding
to the projection of the Weyl vector $w = (0,0,1) \in \Lambda \perp U$ to
$R^\perp\otimes \bbQ$. Then \begin{equation}\label{fh}
h = \Sigma_{R\in \calA\cup\calB} R /3.
\end{equation}
}

\begin{proof} Observe that the right-hand side $h'$ of \eqref{fh}
intersects each curve $R_i$ with multiplicity 1. Since the curves $R_i$
generate $S_X$, $h-h'$ is orthogonal to $S_X$ and hence is equal to zero.
\end{proof}

\subsection{Remark}\label{rk2} Note that the curve $C$ from the proof of 
Lemma \ref{l3} is the cusp curve $C'$ of the quasi-elliptic fibration. In
fact, since it intersects the double components of the reducible fibres, it
is a 2-section. Since the components $F_i$ of fibres and $C'$ generate a
sublattice of finite index in $S_X$,
we must have $C\equiv C'+\sum_i n_iF_i$ for some rational coefficients
$n_i$. Intersecting the both sides with each $F_j$ we get $\sum_i
n_i(F_i\cdot F_j) = 0$. Since the intersection matrix of irreducible
components is semi-definite, this implies that $\sum_i n_iF_i = mF$, where
$F$ is the class of a fibre and $m\in \bbQ$. But now $(C-C')^2 =
-4-2(C\cdot C') = (mF)^2 = 0$ gives $C\sim C'$. Since a smooth rational
curve does not move, we get $C = C'$.

\subsection{Remark}\label{rk3} A quasi-elliptic fibration with 5 fibres of
type $\tilde{D}_4$ has 16 sections
if and only if the Picard lattice has discriminant equal to 4, and hence is
isomorphic to $U\perp D_{20}$.
This follows immediately from the Shioda-Tate formula \cite{Shioda1}. Also,
if $X$ has a quasi-elliptic fibration with five fibres of type
$\tilde{D}_4$ and sixteen sections, then the sixteen sections are
automatically disjoint (note that $m$-torsion sections on an elliptic or a
quasi-elliptic surface over a field of characteristic dividing $m$ are not
necessarily disjoint). Hence the surface contains 42 smooth rational curves
with intersection matrix as described in Lemma \ref{l3}. Let us show this.
Let $C$ be
the cuspidal curve of the fibration and let $F_i =
2E_0(i)+E_1(i)+\ldots+E_4(i), i = 1,\ldots,5,$ be the reducible fibres of
the fibration. The divisor $D = 2C+E_0(1)+\ldots+E_0(4)$ is nef and
satisfies $D^2 = 0$. Thus it defines a genus 1 pencil $|D|$. No sections
intersect $D$, so they are contained in fibres of the fibration. Also the
components $E_1(5),E_2(5),E_3(5),E_4(5)$ do not intersect $D$ and hence are
contained in fibres of $|D|$. We have 16 sections, so there exists
$E_i(5)$, say
$E_1(5)$, such that it intersects at least four of the sections. By
inspection of the list of possible reducible fibres of a genus 1 fibration
we find that $E_1(5)$ together with 4 disjoint sections form a fibre of
type $\tilde{D}_4$.
Now we have 12 remaining sections which intersect one of the components
$E_i(5), i = 2,3,4$. Again we may assume that $E_2(5)$ intersects at least
four of the remaining 12 sections, and hence we find another fibre of type
$\tilde{D}_4$. Continuing in this way we obtain that the 16 sections
together with the curves $E_1(5),E_2(5),E_3(5),E_4(5)$ form four fibres of
type
$\tilde{D}_4$. In particular, they are disjoint.

\subsection{The 168 divisors}\label{168'} Let $r$ be a Leech root as in
Lemma \ref{168}.
The projection $r'$ of $r$ into $S_X \otimes \bbQ$ is a $(-1)$-vector. We
can directly see that each
$r'$ meets exactly 6 members in each family $\calA$ and $\calB$ stated in
Lemma \ref{l3}. For example, if we use the notation as in Remark \ref{rk1}
and take $(\nu_{\Omega} - 4\nu_{4}, 1, 1)$ as $r$, then $r$ meets exactly
twelve Leech roots corresponding to $$E_1, E_5, E_9, E_{10}, E_{11},
4\nu_{\infty} + 4\nu_{0}, L_1, L_2, L_3, L_4, L_5, L_6.$$ We remark that
these twelve curves are mutually disjoint. We set
\begin{equation}\label{class}
l = (2h + \sum_{R \in {\calA}} R )/ 7.
\end{equation}
Then $l \in S_X$, $l^2 = 2$, $l\cdot R = 0$ for any $R \in \calA$, $l\cdot
R = 1$ for any $R \in \calB$ and $$2r' = 2l - (R_1 + \cdot \cdot \cdot +
R_6),$$ where $R_1,..., R_6$ are $(-2)$-curves in $\calA$ which meet $r'$.
Each $r'$ defines an isometry $$s_{r'} : x \to
x + 2(r'\cdot x) r'$$ of $S_X$ which is nothing but the reflection with
respect to the hyperplane perpendicular to $r'$. 

Consider a quasi-elliptic fibration on $X$ with five reducible fibres of
type $\tilde{D}_4$ with 16 sections.
Let $C$ be the cusp curve and let $F$ be the class of a fibre. We now take
as $\calA$ the set of 21 curves
consisting of $C$ and simple components of fibres, and as $\calB$ the set
consisting of 16 sections and multiple
components of fibres. Then we can easily see that \begin{align}\label{168"}
l = C + F,\quad s_{r'}(l) = 5l - 2(R_1+\cdot \cdot \cdot + R_6), \\
s_{r'}(R_i) = 2l - (R_1+\cdot \cdot \cdot + R_6) + R_i.\notag \end{align}

\section{Explicit constructions}

\subsection{Inseparable double covers} Let $\calL$ be a line bundle on a
nonsingular surface $Y$ over a field of
characteristic $2$, and $s\in H^0(Y,\calL^2)$. The pair $(\calL,s)$ defines a
double cover
$\pi:Z\to Y$ which is given by local equations $z^2 = f(x,y)$, where
$(x,y)$ is a system of local parameters and $f(x,y) = 0$ is the local
equation of the divisor of zeroes of the section $s$. Replacing $s$ with
$s+t^2$, where $t\in H^0(Y,\calL)$, we get an isomorphic double cover. The
singular locus of $Z$ is equal to the zero locus of the section $ds\in
H^0(Y,\Omega_Y^1\otimes
\calL^2)$. It is locally given by the common zeroes of the partials of
$f(x,y)$. The canonical sheaf of $X$ is isomorphic to
$\pi^*(\omega_Y\otimes \calL)$. All of these facts are well-known (see, for
example, \cite{CD}, Chap. 0).

We shall consider the special case when $Y =\bbP^2$ and $\calL =
\calO_{\bbP^2}(3)$. The section $s$ is identified with a homogeneous form
$F(x_0,x_1,x_2)$ of degree 6. We assume that $F_6 = 0$ is a reduced plane
curve of degree 6. The expected number $N$ of zeroes of the section $dF_6$
is equal to the second Chern number $c_2(\Omega^1_{\bbP^2}(6))$. The
standard computation gives $c_2 = 21$. So, if $N = 21$, the partial
derivatives of any local equation of $F_6$ at a zero $x$ of $ds$ generate
the maximal ideal $\mathfrak{m}_x$. This easily implies that the cover $Z$
has 21 ordinary double points. So a minimal resolution $X$ of $Z$ has 21
disjoint smooth rational curves.

\subsection{The double plane model}\label{doubleplane} To construct our
surface $X$ as in the Main Theorem \ref{main}, (vi), we take the
inseparable cover
corresponding to the plane sextic $F_6 = 0$ with $dF_6$ vanishing at the 21
points of $\bbP^2$ defined over $\bbF_4$. Consider the sextic defined by
the equation \begin{equation}\label{sextic} F_6 =
x_0x_1x_2(x_0^3+x_1^3+x_2^3) = 0.
\end{equation}
Direct computations show that the partials of $F_6$ vanish at the set
$\calP$ of 21 points of
$\bbP^2(\bbF_4)$:
\[
(1,1,0), (1,0,0),(0,1,0),(0,0,1), (1,1,1), (0,1,1),(1,0,1), \] \[ (1,a,0),
(1,a^2,0), (1,0,a),
(1,0,a^2), , (0,1,a), (0,1,a^2), (1,1,a), \] \[ (1,1,a^2), (1,a,1),
(1,a^2,1), (1,a,a), (1,a,a^2), (1,a^2,1), (1,a^2,a^2), \] where $\bbF_4 =
\{ 0,1,a,a^2 \}$.
Thus the exceptional divisor of a minimal resolution $X$ of the
corresponding double cover $X'\to
\bbP^2$ is a set $\calA$ of 21 disjoint smooth rational curves. Let
$\Sigma(\calP)$ be the blow up of the set $\calP$. Then $X$ is isomorphic
to the normalization of the base change $X'\times_{\bbP^2}\Sigma(\calP)$.
Now let $\check{\calP}$ be the set of 21 lines on $\bbP^2$ defined over
$\bbF_4$. Their equations are $a_0x_0+a_1x_1+a_2x_2 = 0$, where
$(a_0,a_1,a_2)\in \calP$. The 21 lines are divided into three types: three
lines which are components of the sextic $F_6 = 0$; 9 lines which intersect
the sextic at three of its double points; 9 lines which intersect the
sextic at two points of its double points. In the last case each line is a
cuspidal tangent line of the cubic : $x_0^3+x_1^3+x_2^3=0$. The proper
inverse transform $\bar{l}$ of any line in $\Sigma(\calP)$ is a smooth
rational curve with self-intersection $-4$ which is either disjoint from
the proper inverse transform of the sextic, or is tangent to it, or is
contained in it. In each case the pre-image of $\bar{l}$ in $X$ is a
$(-2)$-curve taken with multiplicity 2.
Let $\calB$ be the set of such curves.
Since each line
$\check{p}\in
\check{\calP}$ contains exactly 5 points from $\calP$, and each point $p\in
\calP$ is contained in exactly 5 lines from $\check{\calP}$, the sets
$\calA$ and $\calB$ satisfy property (v) from the Main Theorem. Finally
notice that the map $X\to \bbP^2$ is given by the linear system $|l|$,
where $l$ is defined by \eqref{class}.

\subsection{A switch}\label{switch} A {\it switch} is an automorphism of
$X$ which interchanges the sets $\calA$ and $\calB$. Let us show that it
exists. Consider the pencil of cubic curves generated by the cubics
$x_0x_1x_2 = 0$ and
$x_0^3+x_1^3+x_2^3 = 0$. The set of base points is given by 

\begin{equation}\label{eq1}
(1,1,0), (1,a,0), (1,a^2,0), (1,0,1), (1,0,a), \end{equation} \[
(1,0,a^2), (0,1,1), (0,1,a), (0,1,a^2).
\]

After blowing up the base points, we obtain a rational elliptic surface $V$
with 4
reducible fibres of type $\tilde{A}_2$. The 12 singular points of the four
fibres are the pre-images of the points in $\calP$ not appearing in
\eqref{eq1}.
The surface $X$ is the inseparable cover of the blow-up of $V$ at these
points. It has an elliptic fibration with 4 fibres of type $\tilde{A}_5$.
The base of this fibration is an inseparable double cover of the base of
the elliptic fibration on $V$. Now we find that the Mordell-Weil group is
isomorphic to
$\bbZ/2\bbZ\oplus (\bbZ/3\bbZ)^2$. Nine of these sections are curves from the
set $\calA$ of $(-2)$-curves corresponding to points in $\bbP^2(\bbF_4)$.
They correspond to the base points \eqref{eq1}. Another 9 sections are from
the set $\calB$ of $(-2)$-curves corresponding to lines in
$\bbP^2(\bbF_4)$. They correspond to the lines dual to the points
\eqref{eq1}, i.e. the lines $ax_0+bx_1+cx_2 = 0$, where $(a,b,c)$ is one of
\eqref{eq1}. Fix a zero section
$s_0$ represented by say a $\calA$-curve. Let $F_i = \sum_{k\in
\bbZ/6}E_k(i), i = 1,2,3,4,$ be the reducible fibres such that $
E_k(i)\cdot E_{k+1}(i) = 1$ and $s_0$ intersects $E_0(i)$. It is easy to
see that the 9 sections from $\calA$ intersect the components
$E_0(i),E_2(i),E_4(i)$, three each component. The $\calB$-sections
intersect the components $E_1(i),E_1(i),E_5(i)$, again three each. Also,
the components $E_0(i),E_2(i),E_4(i)$ are $\calB$-curves and the components
$E_1(i),E_1(i),E_5(i)$ are $\calA$-curves. Thus we obtain that
$\calA\cup\calB$ consists of 18 sections and 24 components of fibres. Now
consider the automorphism $T$ of the surface $X$ defined by the translation
by the 2-torsion section $s_1$. Obviously, $s_1$ intersects the components
$E_3(i)$. We see that $T$ interchanges the sets $\calA$ and $B$. Also, note
that if $s_0$
corresponds to $(a_0,a_1,a_2)$, then $s_1$ is a $\calB$-curve corresponding
to the line $a_0x_0 + a_1x_1 + a_2x_2 = 0$. To be more specific, let $s_0$
correspond to the point $p=(1,1,0)$. Then $s_1$ corresponds to the line
$l:x_0+x_1 = 0$. Indeed $l$ is the cuspidal tangent of $F:x_0^3+x_1^3+x_2^3
= 0$ at the point $p$. After we blow-up
$p$, the two sections $s_0$ and $s_1$ intersect at one point of the fibre
represented by $F$. But this could happen only for the 2-section since the
3-torsion sections do not intersect (we are in characteristic 2). Similarly
we see that $T$ sends a section represented by a point $(a,b,c)$ to the
section represented by the dual line $ax_0+bx_1+cx_2 = 0$. Also, it is easy
to see that $T$ transforms a fibre component corresponding to a point to
the component corresponding to the dual line. Thus $T$ is a switch.

We remark that we can directly find a configuration of four fibres of type
$\tilde{A}_5$ and 18 sections on 42 smooth rational curves $\calA \cup
\calB$. This implies the existence of such elliptic fibration (and hence a
switch)
without using the above double plane construction.

\subsection{Mukai's model}\label{mukai} In a private communication to the
second author, S. Mukai suggested that the surface $X$ is isomorphic to a
surface in $\bbP^2\times \bbP^2$ defined by the equations:
\begin{equation}\label{meq}
x_0^2y_0+x_1^2y_1+x_2^2y_2 = 0, \quad x_0y_0^2+x_1y_1^2+x_2y_2^2 = 0.
\end{equation}
To prove this we consider the linear system $|h|$, where $h$ is the divisor
class representing the projection $w'$ of the Weyl vector $w\in L$. Recall
that $h$ is an ample divisor satisfying $h\cdot R = 1$ for any $R\in
\calA\cup \calB$. Consider the quasi-elliptic fibration $|F|$ as in Lemma
\ref{l3}.
We may assume that the set
$\calA$ consists of non-multiple components of fibres and the cusp curve
$C$, and the set $\calB$ consists of 16 sections, and 5 multiple components
of fibres.

Let $D_1=2R_0+R_1+R_2+R_3+R_4$ be a reducible member of $|F|$ and let
$S_1,\ldots, S_4$ be four sections intersecting $R_1$. Consider the
quasi-elliptic pencil $|F'| = |2R_1+S_1+S_2+S_3+S_4|$. The set
$(\calB\setminus \{R_0\})\cup \{C\}$ is the set of irreducible components
of reducible members of $|F'|$. The curve $C' = R_0$ is its cuspidal curve.
Now we check that, for any $R\in \calA\cup\calB$, \[(h-C-F)\cdot R =
(C'+F')\cdot R .\] Since the curves $R$ generate the Picard group of $X$,
we see that \begin{equation}\label{lin}
h = (C+F)+(C'+F').
\end{equation}
The linear system $|C+F|$ defines a degree 2 map $\pi_1:X\to \bbP^2$ which
blows down the curves from the set $\calA$ and maps the curves from the set
$\calB$ to lines. The linear system $|C'+F'|$ defines a degree 2 map
$\pi_2:X\to \bbP^2$ which blows down the curves from the set $\calB$ and
maps the curves from $\calA$ to lines. Let $\phi:X\to \bbP^8$ be the map
defined by the linear system $|h|$.
Using \eqref{lin} one easily sees that $\phi$ maps $X$ isomorphically onto
a surface
contained in the Segre variety $s(\bbP^2\times \bbP^2)$. Let us identify
$X$ with a surface in $\bbP^2\times \bbP^2$. The restriction of the
projections $p_1,p_2:\bbP^2\times\bbP^2\to \bbP^2$ are the maps
$\pi_1,\pi_2$ defined by the linear systems $|C+F|$ and $|C'+F'|$. Let
$l_1,l_2$ be the standard generators of $\Pic(\bbP^2\times \bbP^2)$ and
$[X] = al_1^2+bl_1\cdot l_2+cl_2^2$ be the class of $X$ in the Chow ring of
$\bbP^2\times \bbP^2$. Intersecting $[X]$ with $l_1^2$ and $l_2^2$, we get
$a = c = 2$. Since, $h^2 = 14$ we get
\[14 = [X]\cdot (l_1+l_2)^2 = (2l_1^2+bl_1\cdot l_2+2l_2^2)\cdot
(l_1^2+2l_1\cdot l_2+l_2^2) = 4+2b.\]
Thus
\[[X] = 2l_1^2+5l_1\cdot l_2+2l_2^2 = (2l_1+l_2)\cdot (l_1+2l_2).\] This
shows that $X$ is a complete intersection of two hypersurfaces $V_1,V_2$ of
bidegree $(2,1)$ and $(1,2)$.

The image of each curve from the set $\calA$ (resp. $\calB$) in
$\bbP^2\times\bbP^2$ is a line contained in a fibre of the projection $p_1$
(resp. $p_2$). The fibres of the projection $p_1:V_1\to\bbP^2$ define a
linear system $|L|$ of conics in $\bbP^2$. Let $V$ be the Veronese surface
parametrizing double lines. The intersection $|L|\cap V$ is either subset
of a conic, or is the whole $|L|$. Since $|L|$ contains at least 21 fibres
which are double lines which are not on a conic, we see that $|L|\subset
V$. Thus all fibres of $p_1:V_1\to \bbP^2$ are double lines. This implies
that the equation of $V_1$ can be chosen in the form
\[A_0y_0^2+A_1y_1^2+A_2y_2^2 = 0,\] where the coefficients are linear forms
in $x_0,x_1,x_2$. It is easy to see that the linear forms must be linearly
independent (otherwise $X$ is singular). Thus, after a linear change of the
variables $x_0,x_1,x_2$, we may assume that $A_i = x_i, i= 0,1,2$. 

Now consider a switch $T\in \Aut(X)$ constructed in section \ref{switch}.
Obviously it interchanges the linear systems $|C+F|$ and $|C'+F'|$ and
hence is induced by the automorphisms $s$ of $\bbP^2\times\bbP^2$ which
switches the two factors. This shows that $s(V_1) = V_2$, hence the
equation of $V_2$ is \[y_0x_0^2+y_1x_1^2+y_2x_2^2 = 0.\] 

We remark that the curves
$$(x_0,x_1,x_2) = (a_0,a_1,a_2), \quad a_0y_0^2 + a_1y_1^2 + a_2y_2^2 = 0,
~(a_0,a_1,a_2) \in \bbP^2(\bbF_4)$$ and their images under the switch form
42 smooth rational curves on $X = V_1 \cap V_2$ satisfying property (v) in
the Main theorem.

\subsection{The quartic model}\label{quartic} Consider the quartic curve in
$\bbP^2$ defined by the equation
\begin{equation}\label{dickson}
F_4(x_0,x_1,x_2) =
\end{equation}
\[
x_0^4+x_1^4+x_2^4+x_0^2x_1^2+x_0^2x_2^2+x_1^2x_2^2+x_0x_1x_2(x_0+x_1+x_2) =
0. \]
This is a unique quartic curve defined over $\bbF_2$ which is invariant
with respect to the projective linear group $\PGL(3,\bbF_2) \cong
\PSL(2,\bbF_7)$ (see \cite {Dickson}). Let $Y$ be the quartic surface in
$\bbP^3$ defined by the equation
\begin{equation}\label{qs}
x_3^4+F_4(x_0,x_1,x_2) = 0
\end{equation}
Clearly the group $\PGL(3,\bbF_2)$ acts on $Y$ by projective
transformations leaving the plane $x_3 = 0$ invariant. 

By taking the derivatives, we find that $Y$ has 7 singular points
\begin{equation}\label{sp}
(x_0,x_1,x_2,x_3) = (1,1,0,1),
(1,0,0,1),(0,1,0,1),(0,0,1,1),
\end{equation}
\[
(0,1,1,1),(1,0,1,1),(1,1,1,1).
\]
We denote by $\{ P_i \}_{1\le i \le 7}$ be the set of these singular
points. Each singular point is locally isomorphic to the singular point
$z^4+xy = 0$, i.e. a rational double point of type $A_3$. Let $X$ be a
minimal resolution of $Y$. We claim that $X$ is isomorphic to our surface. 
Note that the points in $\bbP^2$ whose coordinates are the first three
coordinates of singular points \eqref{sp} are the seven points of
$\bbP^2(\bbF_2)$. Let $a_0x_0+a_1x_1+a_2x_2 = 0$ be one of the seven lines
of $\bbP^2(\bbF_2)$. The plane $a_0x_0+a_1x_1+a_2x_2 = 0$ in $\bbP^3$
intersects $Y$ doubly along a conic which passes through 3 singular points
of $Y$.
For example, the plane $x_0+x_1+x_2 = 0$ intersects $Y$ along the conic
given by the equations
\[x_0^2+x_1^2+x_2^2+x_3^2+x_0x_1+x_0x_2+x_1x_2 = x_0+x_1+x_2 = 0.\] We
denote by $\{ C_i' \}_{1\le j \le 7}$ the set of these conics and by $I(j)$
the set of indices $i$ with $P_i \in C_j'$.

Let $R_1^{(i)}+R_2^{(i)}+R_3^{(i)}$
be the exceptional divisor of a singular point $P_i$. We assume that
$R_2^{(i)}$ is the central component. It is easy to check that the proper
inverse transform
$C_j$ of each conic
$C_j'$ in
$X$ intersects $R_2^{(j)}$
with multiplicity 1, if $P_j\in C_i'$, and does not intersect other
components. This easily gives
\begin{equation}\label{quar}
h = 2C_j+\sum_{i\in
I(j)}(R_1^{(i)}+2R_2^{(i)}+R_3^{(i)}),\end{equation} where $h$ is the
pre-image in $X$ of the divisor class of a hyperplane section of $Y$, and
we identify the divisor classes of $(-2)$-curves with the curves.

Observe now that $X$ contains a set $\calA$ of 21 disjoint smooth rational
curves: seven curves $C_i$ and 14 curves $ R_1^{(i)},R_3^{(i)}$. Each of
the seven curves $R_2^{(i)}$ intersects exactly 5 curves from the set
$\calA$
(with multiplicity $1$). We shall exhibit the additional 14 smooth rational
curves which together with these seven curves form a set $\calB$ of 21
disjoint
$(-2)$ curves such that $\calA\cup \calB$ satisfies property (v) of the
Main Theorem.

To do this let us
take a line in
$\bbP^2(\bbF_2)$, for example,
$x_1=0$. Then there
are exactly four points on $\bbP^2(\bbF_2)$ not lying on this line:
$$(x_0,x_1,x_2) =(0,1,0), (1,1,0),(0,1,1),(1,1,1).$$ The plane $H:x_1 + x_3
= 0$ in $\bbP^3$ passes through the 4 singular points of the quartic
surface $Y$:
$$P_1 = (0,1,0,1),\ P_2=(1,1,0,1),\ P_3= (0,1,1,1),\ P_4=(1,1,1,1)$$ and
intersects $Y$ along a quartic curve $Q'$ given by the equations
$$x_0^4+x_2^4+x_0^2x_1^2+x_1^2x_2^2+x_2^2x_0^2 +x_0x_1x_2(x_0+x_1+x_2) =
x_1+x_3=0. $$
It splits into
the union of two conics
$$Q_1': x_0^2+ax_2^2+x_0x_1+ax_1x_2=x_1+x_3=0,$$ $$Q_2':
x_0^2+a^2x_2^2+x_0x_1+a^2x_1x_2=x_1+x_3=0.$$ Let $C'$ be one of the seven
quartic curves $C_i'$'s corresponding to the line $x_1 = 0$. It is given by
the equations
$$x_0^2 + x_2^2 + x_3^2 + x_0x_2 = x_1=0.$$ Then $H$ meets $C'$ at
$q_1=(1,0,a,0),\ q_2=(1,0,a^2,0)$.
Note that $Q_1'$ passes through $P_1,\ldots,P_4, q_1$ and $Q_2$ passes
through $P_1,\ldots,P_4, q_2$.
Each singular point $P_i$ of the quartic is locally isomorphic to $z^4+xy =
0$ and the local equation of the quartic $Q'$ at this point is $z = 0$.
This easily shows that the proper inverse transform of $Q'$ in $X$ consists
of two smooth rational curves $Q_1$ and $Q_2$ each intersects simply the
exceptional divisor at one point lying in different components of
$R_1^{(i)}\cup R_3^{(i)}$. Thus $Q_1$ (or $Q_2$) intersects exactly 5
curves from $\calA$: 4 curves $R_{j}^{(i)}$ ($j = 1$ or $3$, $i =
1,\ldots,4$) and the proper inverse transform of $C'$ in $X$. Also, it is
clear that 14 new quartic curves obtained in this way neither intersect
each other after we resolve the singularities of the quartic nor
intersect the curves $R_2^{(i)}$. This proves the claim.

\subsection{Remark}\label{} The configuration $\calC$ of 14 curves
$C_i,R_2^{(i)}$ is isomorphic to the configuration of points and lines in
$\bbP^2(\bbF_2)$. The group $\PGL(3,\bbF_2)$ acts on the surface $X$ via
its linear action in $\bbP^3$ leaving the hyperplane $x_3 = 0$ fixed. Its
action on the configuration $\calC$ is isomorphic to its natural action on
lines and points.

\section{Automorphisms of $X$}\label{auto} In this section we describe the
group of automorphisms of the surface $X$. First we exhibit some
automorphisms of $X$ and then prove that they generate the group $\Aut(X)$.

\subsection{The group $\PGL(3,\bbF_4)$}\label{group} Consider the double
plane model of $X$ with the branch curve $F_6 = 0$ as in \eqref{sextic}.
The group $G = \PGL(3,\bbF_4)$ acts naturally in the plane
$\bbP^2(\bbF_4)$. For any $g\in G$, let $P_g = F_6(g(x))$. We can see that
$dF(g(x)) = (dF)(g(x))g'(x)$ vanishes at the same set of 21 points on
$\bbP^2(\bbF_4)$. One checks that the linear system of curves of degree 5
vanishing at the 21 points is spanned by the curves \[x_0^4x_1+x_1^4x_0 =
0, \quad x_0^4x_2+x_2^4x_0 = 0,\quad x_2^4x_1+x_1^4x_2 = 0. \] By
integrating we see that $P_g = \alpha_gF_6 + A_g^2$, where $A_g$ is a cubic
form and
$\alpha_g\in k^*$. It is easy to see that the map $\chi:g\to \alpha_g$ is
a homomorphism of groups. The restriction of $\chi$ to  the subgroup
$G' = PSL(3,\bbF_4)$ is trivial, since the latter is a simple group.
Also, it is directly verified that the branch curve does not change under
the transformations given by the diagonal matrices $[a,1,1]$ which
together with $G'$ generate $G$.  Hence, for all $g\in G$, 
$\alpha_g = 1$. Thus the map
$T_g:(z,x)\to (z+A_g(x),g(x))$ is an automorphism of the double plane. It
is easy to verify that
$A_{g'g}(x) = A_{g'}(g(x))+A_{g}(x)$. Since $T_{g'g}(z,x) =
(z+A_{g'g}(x),g'g(x))$ and \[T_{g'}(T_g(z,x)) = T_{g'}(z+A_g(x),g(x)) =
(z+A_{g'}(g(x))+A_{g}(x),g'(g(x)),\] we obtain that the map $g\to T_g$ is a
homomorphism. Since $G$ is simple, it is injective. 

Note that the induced action of $G$ on the sets $\calA$ and $\calB$ is
isomorphic to the action of $G$ on points and lines in $\bbP^2(\bbF_4)$.

\subsection{The 168 Cremona transformations} Let $P_1,\ldots,P_6$ be 6
points in $\bbP^2(\bbF_4)$ such that no three among them are colinear.
Since each smooth conic over $\bbF_4$ contains exactly 5 points (it is
isomorphic to
$\bbP_{\bbF_4}^1$), we see that the set $\calP = \{P_1,\ldots,P_6\}$ is not
on a conic. This allows one to define a unique involutive quintic Cremona
transformation $\Phi$ given by the linear system of curves of degree 5
defined over $\bbF_4$ with double points at the points of $\calP$ (see
\cite{Coo}, Book IV, Chapter VII, \S 4). The transformation $\Phi$ blows
down each conic to a point of $\bbP^2(\bbF_4)$. Let $B:F_6 = 0$ be the
branch curve
\eqref{sextic} of the double cover $X\to \bbP^2$. By adding to $F_6$ the
square of a cubic form, we may assume that each point from $\calP$ is a 
double point of the curve $B$. Since $\Phi^{-1}(\textrm{a line})$ is a
quintic with double points at $\calP$, the image $B' = \Phi(B)$ of $B$ is a
curve of degree
$5\cdot 6-4\cdot 6 = 6$. The image of a conic through 5 points from $\calP$
is a double point of $B'$. Let $F_6'(y_0,y_1,y_2) =0$ be the equation of
$B'$ and let $\Phi$ be given by homogeneous polynomials of degree 5:
\[(y_0,y_1,y_2) =
(f_0(x_0,x_1,x_2),f_1(x_0,x_1,x_2),f_2(x_0,x_1,x_2)).\] Then
\[F_6'(y_0,y_1,y_2) =
F_6(x_0,x_1,x_2)\prod_{i=1}^6q_i(x_0,x_1,x_2),\] where $q_i(x_0,x_1,x_2) =
0$ are the equations of the exceptional conics. Taking the partials we find
\[\sum_{i=0}^2\frac{\partial F_6'}{\partial y_i}\frac{\partial
f_i}{\partial x_j} = \frac{\partial
F_6}{\partial x_j}q+F_6\frac{\partial
q}{\partial x_j}, \quad j = 0,1,2,\]
where $q = \prod_{i=1}^6q_i.$ Let $P\in \bbP^2(\bbF_4)\setminus \calP$. We
know that the partials of $F_6$ vanish at $P$. Adding to $F_6$ a square of
a cubic polynomial we may assume that $F_6(P) = 0$. Since the determinant
of the jacobian matrix $(\frac{\partial f_i}{\partial x_j})$ of $\Phi$ is
invertible outside of the locus $q = 0$, we obtain that the partials of
$F_6'$ vanish at $\Phi(P)$. This shows that the partials of $F_6'$ vanish
at all points of $\bbP^2(\bbF_4)$. Using the same argument as in
\ref{group} we find that $F_6' = \alpha F_6+F_3^2$. Since $\Phi^2$ is the
identity, we get $\alpha^2 = 1$ and hence $\alpha = 1$. Now we can define a
birational transformation of the double plane by the formula
$\tilde{\Phi}(z,x) = (z+F_3,\Phi(x))$. This birational automorphism extends
to a regular automorphism (since $X$ is a minimal model). 

Next let us show that the number of sets $\calP$ is equal to 168. We say
that a subset of $\bbP^2(\bbF_4)$ is
{\it independent} if no three points from it are colinear. Let $N_k$ be the
number of independent subsets of $\bbP^2(\bbF_4)$ of cardinality $k$. We
have
\[N_1
= 21, \quad N_2 = 20N_1/2, \quad N_3 = N_2\cdot 16/3, \quad N_4 = N_3\cdot
9/4,\]
\[ N_5 = N_4\cdot 2/5,\quad
N_6 = N_5/6 = \frac{21\cdot 20\cdot 16\cdot 9\cdot 2}{6!} = 168.\] 

Finally we remark that the above 168 automorphisms act on the Picard
lattice $S_X$ as reflections with respect to 168 $(-4)$-vectors stated in
subsection \ref{168'}. This follows from  equations \eqref{168"}. 

We denote by $N$ the normal subgroup of $\Aut(X)$ generated by 168 involutions. 

\subsection{The automorphism group}\label{Aut}

It is known that the natural map from $\Aut(X)$ to $O(S_X)$ is injective
(\cite{RS2}, \S 8, Proposition 3).
Moreover, $\Aut(X)$ preserves the ample cone, and hence $\Aut(X)$ is a
subgroup of the factor group
$O(S_X)/W(S_X)^{(2)}$, where $W(S_X)^{(2)}$ is the group generated by
$(-2)$-reflections. By the argument in \cite{Kondo}, Lemma 7.3, we can see
that, for any isometry $g$ in $O(S_X)$ preserving an ample class, there
exists an automorphism $\varphi \in N$ such that $g \circ \varphi \in
\Aut(\calD(X))$. This implies that $O(S_X)/W(S_X)^{(2)}$ is a subgroup of a
spilt extension of $N$ by $\Aut(\calD(X))$. Recall that $\Aut(\calD(X))
\cong M_{21} \cdot D_{12}$
(see subsection
\ref{domain}).
Here $M_{21} = \PSL(3,\bbF_4)$.
The Frobenius automorphism of $\bbF_4$ gives an involution on 21 lines and
21 points in $\bbP^2(\bbF_4)$.
This involution induces an isometry $\iota$ of $S_X$ because 42 smooth rational curves generate $S_X$ (Lemma \ref{l4}). The dihedral group $D_{12}$
is generated by $\iota$, a switch and an automorphism of order 3 induced
from a projective transformation of $\bbP^2(\bbF_4)$ given by
$$
\begin{pmatrix}a&0&0\\0&1&0\\0&0&1\end{pmatrix}. $$
We shall show that $\iota$ can not be represented by an automorphism of
$X$. To do this we consider a quasi-elliptic fibration with five fibres of
type $\tilde{D}_4$ from Lemma \ref{l3}. Then $\iota$ fixes three fibres of
type $\tilde{D}_4$ and switches the remaining two fibres. This does not
happen if $\iota$ is realized as an automorphism. 

Thus we conclude that
\[\Aut(X) \cong N\cdot PGL(3,\bbF_4)\cdot 2,\] where the involution 2 is
generated by a switch.

\subsection{Corollary} {\it The finite group $\PGL(3,\bbF_4)\cdot 2$ is
maximal in the following sense.
Let $G$ be a finite group of automorphisms of $X$, then $G$ is conjugate to
a subgroup of
$\PGL(3,\bbF_4)\cdot 2$.}

\begin{proof}
Note that $G$ fixes the vector
$$\tilde{h} = \sum_{g \in G} g^*(h)$$
which is non-zero because $h$ is an ample class and $G$ is an automorphism
group. The vector $\tilde{h}$ is conjugate to a vector in $\calD(X)$ under
$N$ (see subsection \ref{Aut}). This means that $G$ is conjugate to a
subgroup of $\Aut(\calD(X))$. Now the assertion follows. \end{proof} 

\subsection{Conjecture} Let $X$ be a K3 surface over an algebraically
closed field $k$ of characteristic $p$ admitting the group
$\PGL(3,\bbF_4)\cdot 2$ as its group of automorphisms. Then $p = 2$ and $X$
is isomorphic to the surface $X$ from the Main Theorem. 

\section{Proof of the Main Theorem}
\subsection{} By definition, $X$ satisfies property (i) (see subsection
\ref{known}). By Proposition \ref{P1}, $X$ satisfies properties (ii) and
(iii). By Lemma~ \ref{l3}, $X$ satisfies properties (iv) and (v). Property
(vi) from section
\ref{doubleplane}. Property (vii) was proven in \ref{quartic} and property
(viii) in
\ref{mukai}. The group of automorphisms was computed in section \ref{auto}.
The uniqueness follows from the fact that the Artin invariant $\sigma$ is
equal to 1.

It remains to prove the equivalence of properties (i)-(viii). 

(ii) $\Rightarrow$ (i). The components of fibres and the zero section
define the sublattice $U\perp D_{20}\subset \Pic(X)$ of rank 22. If the
assertion is false, $\Pic(X)$ would be unimodular. However, there are no
even unimodular lattices of rank 22 of signature $(1,21)$. 

(i) $\Rightarrow$ (ii) An isotropic vector $f$ from $U$ can be transformed
with the help of $(-2)$-reflections into the class of a fibre of a genus 1
fibration. Let $e_1,\ldots,e_{20}$ be a positive root basis of $D_{20}$. .
Without loss of generality we may assume that $e_1$ is effective. This will
imply that all $e_i$'s are effective. Since $f\cdot e_i = 0$ for all $i$'s,
the irreducible components of $e_i$'s are contained in fibres of the genus
1
fibration. Since the rank of the subgroup of $\Pic(X)$ generated by
irreducible components of fibres is at most 20, we see that all $e_i$'s are
irreducible $(-2)$-curves. They are all contained in one fibre which must
be of type $\tilde{D}_{20}$.
If $f'$ is an isotropic vector from $U$ with $\langle f, f' \rangle = 1$,
the class of $f'-f$ gives the class of a section. A theorem from \S 4 of
\cite{RS2} implies that the fibration is quasi-elliptic. 

(i) $\Leftrightarrow$ (iv) Follows from Remark \ref{rk3}. 

(iv) $\Leftrightarrow$ (v) Follows from Lemma ~\ref{l3}.

(ii) $\Leftrightarrow$ (iii) Follows from Proposition \ref{P1}. 

(vi) $\Leftrightarrow$ (i) Follows from section \ref{doubleplane}.

(vii) $\Leftrightarrow$ (i) Follows from section \ref{quartic}.

(i) $\Leftrightarrow$ (viii) Follows from section \ref{mukai}
and the uniqueness of $X$.

\noindent
This finishes the proof of the Main Theorem.

\subsection{Remark}\label{}
If we take a sublattice $A_2\perp A_2$ in $U\perp \Lambda$, its orthogonal
complement is isomorphic to
$U\perp E_8 \perp E_6 \perp E_6$. This lattice is isometric to the Picard lattice of
the supersingular K3 surface in characteristic 3 with the Artin invariant
1. It is known that this K3 surface is isomorphic to the Fermat quartic
surface
(\cite{Shioda2}, Example 5.2). By using the same method as in this paper we
can see that the projection of the Weyl vector is the class of a hyperplane
section of the Fermat quartic surface, 112 lines on the Fermat quartic
surface can be written in terms of Leech roots and the projective
automorphism group $\text{PGU}(4, \bbF_3)$ of the Fermat quartic surface appears as
a subgroup of
$\Aut(\calD)$.

\subsection{Remark}\label{}
A lattice is called {\it reflective} if its reflection subgroup is of
finite index in the orthogonal group.
The Picard lattice $U\perp D_{20}$ is reflective. This was first pointed
out by Borcherds \cite{Borcherds} and
it is the only known example (up to scaling) of an even reflective lattice of signature
$(1,21)$.

Esselmann \cite{Esselmann} determined the range of possible ranks $r$ of
even reflective lattices of signature $(1,r-1)$: $1 \leq r \leq 20$ or $r =
22$. This is the same as the possible ranks of the Picard lattices of K3
surfaces. Of course, the rank $r = 22$ occurs only when the 
characteristic is positive.


\end{document}